\documentclass[10pt,twoside]{article}
\usepackage{pifont}
\usepackage{bbm}
\usepackage{hyperref}
\usepackage{amsfonts}
\usepackage{amssymb}
\usepackage{amsmath}
\usepackage{amscd}
\usepackage{mathrsfs}
\input amssym.def
\input amssym.tex
\usepackage[all]{xy}
\def\bc{\begin{center}}
\def\ec{\end{center}}
\def\no{\noindent}

\def\c{\cdot}

\def\D{\Delta}

\def\o{\otimes}

\def\r{\rho}

\def\v{\varepsilon}

\begin{document}
\begin{center}
{\Large\bf Rota-Bater paried comodule and Rota-Bater paired Hopf module$^{\ast}$}

\vspace{4mm}

{ Huihui Zheng, Yuxin Zhang, Liangyun Zhang$^{\ast\ast}$}
$$\begin{array}{rllr}
&$~College of Science, Nanjing Agricultural University,
Nanjing 210095, China$
\end{array}$$
\begin{figure}[b]
\rule[-2.5truemm]{5cm}{0.1truemm}\\[2mm]
\no $^{\ast}$This work is supported by Natural Science Foundation (11571173).

$^{\ast\ast}$Corresponding author: zlyun@njau.edu.cn
\end{figure}
\begin{minipage}{144.6mm}

\vskip1cm \footnotesize{\textbf{Abstract}:
In this paper, we introduce the conception of Rota-Baxter paired comodules, which is dual to Rota-Baxter paired modules in \cite{Zheng}. We mainly discuss some properties of Rota-Baxter paired comodules, especially we give the characterization of generic Rota-Baxter paired comodules, which has important application for the construction of Rota-Baxter comodules. Moreover, we construct Rota-Baxter paired comodules on Hopf algebras, weak Hopf algebras, weak Hopf modules, dimodules, relative Hopf modules and Rota-Baxter paired comodules. And then we finally introduce the conception of Rota-Baxter paired Hopf modules by combining Rota-Baxter paired module with Rota-Baxter paired comodule,  and give the structure theorem of generic Rota-Baxter paired Hopf modules.
\vspace{3mm}}\vspace{3mm}\\
\textbf{2020 AMS Classification Number}: 16T05; 16T15; 17B38.
\end{minipage}
\end{center}

\vspace{4mm}

\vskip1mm

\begin{center}
{\bf \S1\quad Introduction}
\end{center}

A Rota-Baxter algebra (first known as a Baxter algebra) is an algebra $A$ with a linear operator $P$ on $A$ that satisfies the Rota-Baxter identity
\begin{equation*}
P(x)P(y)=P(P(x)y)+P(xP(y))+\lambda P(xy),\ \text{for all}\ x,y\in A,
\end{equation*}
where $\lambda\in k$ (the field), called the weight~\cite{Baxter,Rota}.

Rota-Baxter algebra originated from the 1960 paper \cite{Baxter} of Baxter based on his probability study to understand Spitzer's identity in fluctuation theory. It wasn't long before the concept attracted the attention of many mathematicians, especially Rota, whose fundamental papers around 1970 brought the subject into the areas of algebra and combinatorics. In \cite{Aguiar}, a connection with mathematical physics was also established that related a Rota-Baxter algebra of weight 0 to the associative analog of classical Yang-Baxter equation.

To study the representations of Rota-Baxter algebras, the authors in \cite{Guo} introduced the conception of Rota-Baxter modules related to the ring of Rota-Baxter operators. By the definition, a Rota-Baxter module over a Rota-Baxter algebra $(A,P)$ is a pair $(M,T)$ where $M$ is a (left) $A$-module and $T:M\to M$ a $k$-linear operator such that
\begin{equation*}
P(a)\cdot T(m)=T(P(a)\cdot m)+T(a\cdot T(m))+\lambda T(a \cdot m), \ \ \text{for all}~a\in A, m\in M.
\end{equation*}

Later, Rota-Baxter paired modules were introduced in \cite{Zheng}, without requiring $(A,P)$ to be a Rota-Baxter algebra, which is a natural generalization of Rota-Baxter modules. Many properties of Rota-Baxter modules, even of Rota-Baxter algebras, are naturally generalized to Rota-Baxter paired modules. Rota-Baxter paired module has broader connections and applications, especially to Hopf algebras. We have constructed a large number of Rota-Baxter paired modules from Hopf algebra related structures in \cite{Zheng}.

Representation theory of coalgebras and comodules is a very extensive subject. On the basis of the comodule theory, we can naturally consider Rota-Baxter operators on comodules. In this paper, we naturally introduce the conception of Rota-Baxter paired comodules, which is dual to Rota-Baxter paired modules, and give some properties of Rota-Baxter paired comodules. In addition, we also give its construction from Hopf algebra related coalgebras and comodules.

A Hopf module on a bialgebra $H$ is also an $H$-module and an $H$-comodule, whose action and coaction satisfy a compatibility condition.

As is well known, the structure theorem on Hopf modules is concerned by many experts and scholars. Especially, this structure theorem can describe the integrals of Hopf algebras. Combining Rota-Baxter paired modules and Rota-Baxter paired comodules, we can naturally introduce the conception of Rota-Baxter paired Hopf modules, and study Rota-Baxter operator on it, and prove its structure theorem.

 This article is organized as follows. In Section 2, we recall the definition of Rota-Baxter coalgebras, and then give the notion
of Rota-Baxter paired comodules, which is dual to Rota-Baxter paired modules in \cite{Zheng}. Moreover, we provide a large number of examples of Rota-Baxter paired comodules. In Section 3, we discuss some properties of Rota-Baxter paired comodules, especially we give the characterization of generic Rota-Baxter paired comodules (see Theorem 3.1), which has an important application for the construction of Rota-Baxter comodules. In Section 4, we construct Rota-Baxter paired comodules on Hopf algebras, weak Hopf algebras, weak Hopf modules, dimodules, relative Hopf modules and Rota-Baxter paired comodules, respectively. Especially, we find some Rota-Baxter coalgebras and Rota-Baxter paired comodules by applying (co)integral in bialgebra, antipode and idempotent element in (weak) Hopf algebras and $R$-matrix in quasitriangular Hopf algebra. In Section 5, we construct pre-Lie comodules from Rota-Baxter paired comodules. In Section 6, we introduce the conception of Rota-Baxter paired Hopf modules by combining Rota-Baxter paired module with Rota-Baxter paired comodule, and give the structure theorem of generic Rota-Baxter Hopf modules.

 Throughout this paper, let $k$ be a fixed field. Unless otherwise specified, linearity, modules and $\otimes$ are all meant over $k$. And we freely use the Hopf algebras terminology introduced in \cite{Sweedler}. For a coalgebra $C$, we write its comultiplication $\Delta(c)$ with $c_{1}\otimes c_{2}$, for any $c\in C$; for a left $C$-comodule $M$, we denote its coaction by $\rho(m)=m_{(-1)}\otimes m_{(0)}$, for any $m\in M$; for a right $C$-comodule $M$, we denote its coaction by $\rho(m)=m_{(0)}\otimes m_{(1)}$, for any $m\in M$, in which we omit the summation symbols for convenience.

\vspace{3mm}

\begin{center}
{\large\bf \S2\quad Rota-Baxter coalgebras and Rota-Baxter paired comodules}
\end{center}

In this section, we firstly recall the definition of Rota-Baxter coalgebras, and then give the notion of Rota-Baxter paired comodules, which is dual to the definition of Rota-Baxter paired modules in \cite{Zheng}. Moreover, we provide a large number of examples of Rota-Baxter paired comodules.

 \vspace{2mm}

{\bf 2.1\quad Rota-Baxter coalgebras}

 \vspace{2mm}

 \textbf{Definition 2.1} Let $(C,\Delta,\varepsilon)$ be a coalgebra. We call $(C, P)$ a {\bf Rota-Baxter coalgebra of weight $\lambda$} \cite{Jian}, if the linear map $P:C\rightarrow C$ satisfies the following
$$
(P\otimes P)\Delta=(P\otimes id)\Delta P+(id\otimes P)\Delta P+\lambda \Delta P,
$$
where $\lambda\in k$ and $id$ denotes the identical map.
\vspace{4mm}

We refer the reader to \cite{Jian} for further discussions on Rota-Baxter coalgebra and only give the following simple
examples of Rota-Baxter coalgebras which will be revisited later.

 \vspace{2mm}

 \textbf{Example 2.2} \ (a) Let $C$ be an augmented coalgebra, that is, there exists
a coalgebra homomorphism $f: k\rightarrow C$. Then, it is easy to see that $f(1_k)$ is a group-like element in $C$. So,
$(C, P)$ is a Rota-Baxter coalgebra of weight $-1$.

Here $P$ is given by
$$
P: C\rightarrow C, c\mapsto\varepsilon(c)f(1_k).
$$

Furthermore, if $C$ is a bialgebra with the unit $\mu$, then, $C$ is an augmented coalgebra since $\mu: k\rightarrow C$ is a coalgebra map. So, $(C, P)$ is a Rota-Baxter coalgebra of weight $-1$ , where
$$
P: C\rightarrow C, c\mapsto\varepsilon(c)1_C.
$$

(b) Let $C$ be a coalgebra and $\lambda\in C^\ast$ (the linear dual space of $C$). Define
 $$
 P: C\rightarrow C, c\mapsto \lambda(c_1)c_2.
 $$

Then, $(C, P)$ is a Rota-Baxter coalgebra of weight $-1$ if and only if $\lambda^2=\lambda$, that is, $\lambda(c_1)\lambda(c_2)=\lambda(c)$ for all $c\in C$.

 \vspace{2mm}

{\bf 2.2\quad Rota-Baxter  paired comodules}

 \vspace{2mm}

 \textbf{Definition 2.3} Let $(C,\Delta,\varepsilon)$ be a coalgebra, and $M$ a left $C$-comodule with coaction $\rho$. A pair $(P,T)$ of linear maps $P:C\rightarrow C$ and $T:M\rightarrow M$ is called a {\bf Rota-Baxter paired operator of weight $\lambda$ on $M$} if
\begin{equation*}
(P\otimes T)\rho=(P\otimes id)\rho T+(id\otimes T)\rho T+\lambda \rho T.
\end{equation*}

We also call the triple $(M,P,T)$ a {\bf Rota-Baxter paired (left) $C$-comodule of weight $\lambda$}. Given a linear map $T:M\rightarrow M$, if $(M,P,T)$ is a Rota-Baxter paired $C$-comodule of weight $\lambda$ for every linear map $P:C\rightarrow C$, then $(M,T)$ is called a {\bf generic Rota-Baxter paired $C$-comodule of weight $\lambda$}.

\vspace*{2mm}

\textbf{Example 2.4} (1) Let $C$ be a coalgebra, regarded also as a left $C$-comodule via its comultiplication $\Delta$. If $(C,P)$ is a Rota-Baxter coalgebra of weight $\lambda$, then
 $(C,P,P)$ is a Rota-Baxter paired $C$-comodule of weight $\lambda$.

 In particular, if $C$ is a bialgebra with the unit $\mu$, then, by Example 2.2,  $(C,P,P)$ is  a Rota-Baxter paired $C$-comodule of weight $\lambda$, where $
P: C\rightarrow C, c\mapsto\varepsilon(c)1_C.$

(2) Let $(M,P,T)$ be a Rota-Baxter paired $C$-comodule of weight $\lambda$. Then, for any $\mu\in k$, $(M,\mu P, \mu T)$ is a Rota-Baxter paired $C$-comodule of weight $\lambda\mu$.

(3) Let $H$ be a bialgebra, and $(M,P,T)$ a Rota-Baxter paired $H$-comodule of weight $\lambda$. If $P$ is idempotent and a bialgebra homomorphism from $H$ to $H$, then $(H\o M, \r, P, T^{'})$ is a Rota-Baxter paired $H$-comodule of weight $\lambda$, where $T^{'}:H\o M\rightarrow H\o M$ and $\r: H\o M\rightarrow H\o H\o M$ are defined by $T^{'}(h\o m)=P(h)\o T(m)$ and $\r(h\o m)=h_1m_{(-1)}\o h_2\o m_{(0)}$, respectively.

In fact, it is easy to prove that $(H\o M, \r)$ is a left $H$-comodule. Moreover, for any $h\in H, m\in M,$ we have
\begin{eqnarray*}
&&P((h\o m)_{(-1)})\o T'((h\o m)_{(0)})\\
&&=P(h_1m_{(-1)})\o P(h_2)\o T(m_{(0)})\\
&&=P(h_1)P(m_{(-1)})\o P(h_2)\o T(m_{(0)})\\
&&=P(h_1)P(T(m)_{(-1)})\o P(h_2)\o T(m)_{(0)}+P(h_1)T(m)_{(-1)}\o P(h_2)\o T(T(m)_{(0)})\\
&&+\lambda P(h_1)T(m)_{(-1)}\o P(h_2)\o T(m)_{(0)},\\
&&P(T'(h\o m)_{(-1)})\o T'(h\o m)_{(0)}+T'(h\o m)_{(-1)}\o T'(T'(h\o m)_{(0)})\\
&&+\lambda T'(h\o m)_{(-1)}\o T'(h\o m)_{(0)}\\
&&=P(P(h)_1 T(m)_{(-1)})\o P(h)_2\o T(m)_{(0)}+P(h)_1 T(m)_{(-1)}\o P(P(h)_2)\o T(T(m)_{(0)})\\
&&+\lambda P(h)_1 T(m)_{(-1)}\o P(h)_2\o T(m)_{(0)}\\
&&=P(P(h)_1) P(T(m)_{(-1)})\o P(h)_2\o T(m)_{(0)}+P(h)_1 T(m)_{(-1)}\o P(P(h)_2)\o T(T(m)_{(0)})\\
&&+\lambda P(h)_1 T(m)_{(-1)}\o P(h)_2\o T(m)_{(0)}\\
&&=P(h_1)P(T(m)_{(-1)})\o P(h_2)\o T(m)_{(0)}+P(h_1)T(m)_{(-1)}\o P(h_2)\o T(T(m)_{(0)})\\
&&+\lambda P(h_1)T(m)_{(-1)}\o P(h_2)\o T(m)_{(0)},
\end{eqnarray*}
so, by Definition 2.3, we know that $(H\o M, \r, P, T^{'})$ is a Rota-Baxter paired $H$-comodule of weight $\lambda$.

(4) Let $M$ be a left $C$-comodule with the coaction $\rho$, and $V$ a vector space. Then, $M\otimes V$ has a left $C$-comodule structure, whose comodule structure map is given by $\rho\otimes id$. So, if $(M, P, T)$ is a Rota-Baxter paired $C$-comodule of weight $\lambda$, we easily see that
 $(M\otimes V, P, T\otimes id)$ is a Rota-Baxter paired $C$-comodule of weight $\lambda$.

 In particular, if $(C,P)$ is a Rota-Baxter coalgebra of weight $\lambda$, then $(C\o V, P, P\otimes id)$ is a Rota-Baxter paired $C$-comodule of weight $\lambda$.

Furthermore, if $(M, T)$ is a generic Rota-Baxter paired $C$-comodule of weight $\lambda$,  $(M\otimes V, T\otimes id)$ is also a generic Rota-Baxter paired $C$-comodule of weight $\lambda$.

(5) Let $M$ be a left $C$-comodule, and $T$ an idempotent epimorphism in End$(M)$. Then, $(M,id,T)$ is a Rota-Baxter paired $C$-comodule of weight $-1$.

In fact, for any $m\in M,$ we have
\begin{eqnarray*}
((id\otimes id)\r T^2+(id\otimes T)\r T^2-\r T^2)(m)&=&((id\otimes id)\r T+(id\otimes T)\r T-\r T)(m)\\
&=&(id\otimes T)\r T(m).
\end{eqnarray*}

Since $T$ is an idempotent epimorphism, we have
$$(id\otimes T)\r=(id\otimes id)\r T+(id\otimes T)\r T-\r T.
$$

Hence, $(M,id,T)$ is a Rota-Baxter paired $C$-comodule of weight $-1$.

\vspace*{2mm}

A {\it Rota-Baxter paired $C$-subcomodule $N$} of a Rota-Baxter paired $C$-comodule $(M, P, T)$ is a
$C$-subcomodule of $M$ such that $T(N)\subseteq N$. A {\it Rota-Baxter paired comodule map} $f: (M, P, T) \rightarrow (M^\prime, P^\prime, T^\prime) $ of the same weight $\lambda$ is a $C$-comodule map
such that $fT=T^{\prime}f$.

\vspace*{2mm}

\textbf{Proposition 2.5}\ \ Let $f : (M, P, T) \rightarrow (M^\prime, P^\prime, T^\prime)$ be a Rota-Baxter paired comodule map of weight $\lambda$. Then the following conclusions hold.

(a) Ker$f$ is a Rota-Baxter paired $C$-subcomodule of $M$.

(b) If $K$ is a Rota-Baxter paired $C$-subcomodule of $M$, then $f(K)$ is a Rota-Baxter paired $C$-subcomodule of $M^\prime$.

In particular, if $T$ is $C$-colinear, then $T(M)$ is a Rota-Baxter paired $C$-subcomodule of $M$.

(c) If $L$ is a Rota-Baxter paired $C$-subcomodule of $M^\prime$, then $f^{-1}(L)$ is a Rota-Baxter paired $C$-subcomodule of $M$.

{\bf Proof.}
(a) Since $f$ is a $C$-comodule map, Ker$f$ is a $C$-subcomodule of $M^\prime$. In addition, for any $x\in$Ker$f$, $fT(x)=T^\prime f(x)=0$, so $T($Ker$f)\subseteq $Ker$f$. Hence Ker$f$ is a Rota-Baxter paired $C$-subcomodule of $M$.

(b) It is obvious that $f(K)$ is a subcomodule of $M^\prime$, so, we have only to verify that $T^\prime(f(K))\subseteq f(K)$.
 Since $T(K)\subseteq K$, and $fT=T^\prime f$, we have
$T^\prime f(K)=fT(K)\subseteq f(K)$.

(c) We consider the composition $\pi f$ of comodule maps $\pi$ and $f$, where $\pi: N\rightarrow N/L$ is a projection. By (a), we know that Ker$f$ is a $C$-subcomodule of $M$, so, Ker$(\pi f)=f^{-1}(L)$ and a subcomodule of $M$.
In addition,  we have
$fT(f^{-1}(L))=T^\prime f(f^{-1}(L))=T^\prime(L)\subseteq L$.
So, we can get $T(f^{-1}(L))\subseteq f^{-1}(L)$.
\hfill $\square$

\vspace*{2mm}

\begin{center}
{\bf \S3\quad Some properties of Rota-Baxter paired comodules}
\end{center}
\vspace*{2mm}

In this section, we will discuss some properties of Rota-Baxter paired comodules.

Recall that a linear operator $T:M\rightarrow M$ is called quasi-idempotent \cite{Zheng} of weight $\lambda$ if $T^2=-\lambda T$. We have the following characterization of generic Rota-Baxter paired comodules, which has important application for the construction of Rota-Baxter comodules.

\textbf{Theorem 3.1}\ \ Let $C$ be a coalgebra, and $M$ a left $C$-comodule. If there exists a colinear map $T: M\rightarrow M$, then the following are equivalent.

(1) $(M,T)$ is a generic Rota-Baxter paired $C$-comodule of weight $\lambda$.

(2) There is a linear operator $P: C\rightarrow C$ such that $(M,P,T)$ is a Rota-Baxter paired $C$-comodule of weight $\lambda$.

(3) $T$ is quasi-idempotent of weight $\lambda$.

{\bf Proof.} Under the $C$-colinearity condition of $T$, for any linear operator $P:A\rightarrow A$ and $m\in M$, we have

$P(m_{(-1)})\otimes T(m_{(0)})=P(T(m)_{(-1)})\otimes T(m)_{(0)}+T(m)_{(-1)}\o T(T(m)_{(0)})+\lambda T(m)_{(-1)}\o T(m)_{(0)}$

$\Longleftrightarrow$
$P(m_{(-1)})\otimes T(m_{(0)})=P(m_{(-1)})\otimes T(m_{(0)})+m_{(-1)}\o T^2(m_{(0)})+\lambda m_{(-1)}\o T(m_{(0)})$

$\Longleftrightarrow$
$0=m_{(-1)}\o T^2(m_{(0)})+\lambda m_{(-1)}\o T(m_{(0)})$.

If (1) holds, applying $\v\o id$ to both sides of the above equation, we get $T^2=-\lambda T$. Conversely, if $T^2=-\lambda T$, it is obvious that (1) holds.

In a similar way, we can prove that $(2)\Longleftrightarrow (3).$
\hfill $\square$

\vspace*{2mm}

\textbf{Proposition 3.2}~ Let $M$ be a left $C$-comodule. Then, there exists a left $C$-comodule map $T: M\rightarrow M$ such that $(M,T)$ is a generic Rota-Baxter $C$-comodule of weight $-1$, if and only if there is a $C$-comodule direct sum decomposition $M=M_1\oplus M_2$ such that $T: M\rightarrow M_1\subseteq M$ is the project of $M$ onto $M_1$: $T(m_1+m_2)=m_1$ for $m_1\in M_1$ and $m_2\in M_2$.

\textbf{Proof.} Suppose $M$ has a direct sum decomposition $M=M_1\oplus M_2$ of $C$-comodules, where $M_1$ and $M_2$ are subcomodule of $M$. Then the projection $T$ of $M$ onto $M_1$ is idempotent, since, for $m=m_1+m_2\in M$ with $m_1\in M_1$ and $m_2\in M_2$, we have
$T^2(m)=T^2(m_1+m_2)=T(m_1)=m_1=T(m)$.

Furthermore, we have
\begin{eqnarray*}
(id\o T)\r(m)&=&(id\o T)\r(m_1+m_2)\\
&=&m_{1(-1)}\o T(m_{1(0)}+0)+m_{2(-1)}\o T(0+m_{2(0)})\\
&=&m_{1(-1)}\o m_{1(0)}=\r(m_1)\\
&=&\r T(m),
\end{eqnarray*}
so, $T$ is a left $C$-comodule map. Again by Theorem 3.1, we know $(M,T)$ is a generic Rota-Baxter paired $C$-comodule of weight $-1$.

Conversely, if $(M, T)$ is a generic Rota-Baxter paired $C$-comodule of weight $-1$ and $T$ a left $C$-comodule map, then by Theorem 3.1, we know $T$ is idempotent.

Let $M_1=T(M)$ and $M_2=(id-T)(M)$. Because $T$ is a left $C$-comodule map, both $M_1$ and $M_2$ are subcomodule of $M$. Also, for any $m\in M$, $m=T(m)+(id-T)(m)$, so $M=M_1+M_2$. Furthermore, if $n\in M_1\cap M_2$, then $n=T(x)=(id-T)(y)$, for some $x,y\in M$. Thus $n=T(x)=T^2(x)=T(id-T)(y)=(T-T^2)(y)=0$. Therefore $M=M_1\oplus M_2$.

Finally, since $m=T(m)+(id-T)(m)$ is the decomposition of $m\in M$ as $m=m_1+m_2$ with $m_1\in M_1$ and $m_2\in M_2$, we see that $T$ is the projection of $M$ onto $M_1$.
\hfill $\square$

\vspace{2mm}

 \textbf{Proposition 3.3} Let $M$ be a $C$-comodule and $P : C\rightarrow C$, $T : M\rightarrow M$ linear maps. Then $(M,P,T)$ is a Rota-Baxter paired $C$-comodule of weight $\lambda \neq 0$ if and only if there is a map $f:M \rightarrow C\otimes M$ such that
 \begin{equation*}
(P\otimes T)\rho =fT, \quad  (\overline{P}\otimes \overline{T})\rho =-f\overline{T},
\end{equation*}
where $\overline{P}=-P-\lambda id$ and $\overline{T}=-T-\lambda id.$

\vspace{2mm}

\textbf{Proof.} Let $(M,P,T)$ be a Rota-Baxter paired $C$-comodule of weight $\lambda$. Then, we have
\begin{equation}
(P\otimes T)\rho = (P\otimes id)\rho T+(id\otimes T)\rho T+\lambda \rho T.  \nonumber
\end{equation}

Let $f=(P\otimes id)\rho +(id\otimes T)\rho +\lambda \rho $. Then the above equation gives
$$
(P\otimes T)\rho=fT,
$$
so, we obtain
\begin{equation*}
(\overline{P}\otimes \overline{T})\rho =-f\overline{T}.
\end{equation*}

Now we consider the converse. Suppose that there exists a map $f:M \rightarrow C\otimes M$  such that
$(P\otimes T)\rho =fT$ and $(\overline{P}\otimes \overline{T})\rho =-f\overline{T}$. Then we have
\begin{eqnarray*}
-\lambda f&=&f \overline{T}+f T=(P\otimes T)\rho-(\overline{P}\otimes \overline{T})\rho \\
&=&(P\otimes T)\rho-((-\lambda id-P)\otimes (-\lambda id- T))\rho\\
&=&(P\otimes T)\rho-(\lambda^2 id\otimes id+\lambda id \otimes T+\lambda P\otimes id+P\otimes T)\rho\\
&=&-\lambda(\lambda id\otimes id+id\otimes T+P\otimes id)\rho.
\end{eqnarray*}

So, we get
$$
f=(P\otimes id)\rho +(id\otimes T)\rho +\lambda \rho.
$$

Furthermore, $(P\otimes T)\rho=((P\otimes id)\rho +(id\otimes T)\rho +\lambda \rho) T$. Hence $(M,P,T)$ is a Rota-Baxter paired $C$-comodule of weight $\lambda\neq 0$.
\hfill $\square$

\vspace*{2mm}

By the above definition of $\overline{T}$ and $\overline{P}$ in Proposition 3.3, there are also the following relationships.

\vspace*{2mm}

\textbf{Proposition 3.4}\ \  Let $(M,P,T,\rho)$ be a Rota-Baxter paired $C$-comodule of weight $\lambda$. Then we have
\begin{eqnarray*}
&&P(m_{(-1)})\o \overline{T}(m_{(0)})=T(m)_{(-1)}\o \overline{T}(T(m)_{(0)})+P(\overline{T}(m)_{(-1)})\o \overline{T}(m)_{(0)},\\
&&\overline{P}(m_{(-1)})\o T(m_{(0)})=\overline{P}(T(m)_{(-1)})\o T(m)_{(0)}+\overline{T}(m)_{(-1)}\o T(\overline{T}(m)_{(0)})
\end{eqnarray*}
for any $m\in M$.

{\bf Proof.} By using the compatible condition of Rota-Baxter paired comodules, we can directly verify that the equation holds.
\hfill $\square$

\vspace*{2mm}

The following result shows that how close it is for an idempotent Rota-Baxter operator of comodule to have weight $-1$.

\vspace*{2mm}

\textbf{Proposition 3.5}\ \ Let $(M,P,T)$ be a Rota-Baxter paired $C$-comodule of weight $\lambda$.

(a) If $T$ is idempotent, then $(1+\lambda)T(m)_{(-1)}\o T(T(m)_{(0)})=0$, for any $m\in M$.

(b)  If $P$ and $T$ are idempotent, then $(1+\lambda)P(T(m)_{(-1)})\o T(m)_{(0)}=0$, for any $m\in M$.

(c)  If $P$ and $T$ are idempotent, then $(1+\lambda)(P(T(m)_{(-1)})\o T(m)_{(0)}-\lambda T(m)_{(-1)}\o T(m)_{(0)})=0$, for any $m\in M$.

 Thus an idempotent Rota-Baxter operator of comodule must have weight $-1$.

{\bf Proof.} (a) Since $T^2=T$, for any $m\in M$ we obtain
\begin{eqnarray*}
&&P(m_{(-1)})\o T(m_{(0)})=P(m_{(-1)})\o T(T(m_{(0)}))\\
&=&P(T(m)_{(-1)})\o T(T(m)_{(0)})+T(m)_{(-1)}\o T^2(T(m)_{(0)})+\lambda T(m)_{(-1)}\o T(T(m)_{(0)})\\
&=&P(T(m)_{(-1)})\o T(T(m)_{(0)})+T(m)_{(-1)}\o T(T(m)_{(0)})+\lambda T(m)_{(-1)}\o T(T(m)_{(0)}),
\end{eqnarray*}
that is, $(P\o T)\rho=(P\o T)\rho T+(id\o T)\rho T+\lambda(id\o T)\rho T$.

Applying $T^2=T$ to the above equality, we have $(P\o T)\rho T=(P\o T)\rho $. Thus $(1+\lambda)T(m)_{(-1)}\o T(T(m)_{(0)})=0$ for any $m\in M$.

(b) Since $P^2=P$, for any $m\in M$ we obtain
\begin{eqnarray*}
&&P(m_{(-1)})\o T(m_{(0)})=P(P(m_{(-1)}))\o T(m_{(0)})\\
&=&P^2(T(m)_{(-1)})\o T(m)_{(0)}+P(T(m)_{(-1)})\o T(T(m)_{(0)})+\lambda P(T(m)_{(-1)})\o T(m)_{(0)}\\
&=&P(T(m)_{(-1)})\o T(m)_{(0)}+P(T(m)_{(-1)})\o T(T(m)_{(0)})+\lambda P(T(m)_{(-1)})\o T(m)_{(0)},
\end{eqnarray*}
that is, $(P\o T)\rho=(P\o id)\rho T+(P\o T)\rho T+\lambda(P\o id)\rho T$.

Applying $T^2=T$ to the above equality, we have $(P\o T)\rho T=(P\o T)\rho $.
Thus $(1+\lambda)P(T(m)_{(-1)})\o T(m)_{(0)}=0$ for any $m\in M$.

(c) By (a) and (b), we can prove that (c) holds.
\hfill $\square$

\begin{center}
{\bf \S4\quad Constructios of Rota-Baxter paired comodule}
\end{center}

In this section, we construct Rota-Baxter paired comodule by deformation, direct sum, Hopf algebra, Rota-Baxter paired module, Hopf module, co-Hopf module, dimodule.

Firstly, we construct Rota-Baxter paired comodule by deformation.

\vspace*{2mm}

{\bf 4.1\quad The construction on Hopf algebras}

 \vspace{2mm}

\textbf{Proposition 4.1} Let $C$ be a coalgebra, and $M$ a left $C$-comodule. Define two maps
$T: M\rightarrow M$ and $P: C\rightarrow C$ by $T(m)=\chi (m_{(-1)})m_{(0)}$ and $P(c)=\chi(c_1)c_2$, for any $m\in M, c\in C$, respectively. Then, $(M,P,T)$ is a Rota-Baxter paired $C$-comodule of weight $-1$ if $\chi\in C^\ast$ is idempotent under the convolution product.

\textbf{Proof.} For any $m\in M$, we have
\begin{eqnarray*}
P(m_{(-1)})\o T(m_{(0)})&=&\chi(m_{(-1)1})m_{(-1)2}\o \chi(m_{(0)(-1)})m_{(0)(0)}\\
&=&\chi(m_{(-1)1})m_{(-1)2}\o \chi(m_{(-1)3})m_{(0)}.
\end{eqnarray*}

Moreover, for any $m\in M$, $\rho(T(m))=\rho(\chi (m_{-1})m_{0})=\chi (m_{-1})m_{(0)(-1)}\o m_{(0)(0)}=\chi (m_{(-1)1})m_{(-1)2}\o m_{(0)}$, so, we obtain that
\begin{eqnarray*}
&&P(T(m)_{(-1)})\o T(m)_{(0)}+T(m)_{(-1)}\o T(T(m)_{(0)})-T(m)_{(-1)}\o T(m)_{(0)}\\
&&=\chi (m_{(-1)1})P(m_{(-1)2})\o m_{(0)}+\chi (m_{(-1)1})m_{(-1)2}\o T(m_{(0)})-\chi (m_{(-1)1})m_{(-1)2}\o m_{(0)}\\
&&=\chi (m_{(-1)1})\chi(m_{(-1)2})m_{(-1)3}\o m_{(0)}+\chi (m_{(-1)1})m_{(-1)2}\o \chi (m_{(0)(-1)})m_{(0)(0)}-\chi (m_{(-1)1})m_{(-1)2}\o m_{(0)}\\
&&=\chi (m_{(-1)1})\chi(m_{(-1)2})m_{(-1)3}\o m_{(0)}+\chi(m_{(-1)1})m_{(-1)2}\o \chi(m_{(-1)3})m_{(0)}-\chi (m_{(-1)1})m_{(-1)2}\o m_{(0)}\\
&&=\chi^2(m_{(-1)1})m_{(-1)2}\o m_{(0)}+\chi(m_{(-1)1})m_{(-1)2}\o \chi(m_{(-1)3})m_{(0)}-\chi (m_{(-1)1})m_{(-1)2}\o m_{(0)}\\
&&=\chi(m_{(-1)1})m_{(-1)2}\o m_{(0)}+\chi(m_{(-1)1})m_{(-1)2}\o \chi(m_{(-1)3})m_{(0)}-\chi (m_{(-1)1})m_{(-1)2}\o m_{(0)}\\
&&=\chi(m_{(-1)1})m_{(-1)2}\o \chi(m_{(-1)3})m_{(0)}.
\end{eqnarray*}

Hence $P(m_{(-1)})\o T(m_{(0)})=P(T(m)_{(-1)})\o T(m)_{(0)}+T(m)_{(-1)}\o T(T(m)_{(0)})-T(m)_{(-1)}\o T(m)_{(0)}$, $(M,P,T)$ is a Rota-Baxter paired $H$-comodule of weight $-1$.
\hfill $\square$

\vspace*{2mm}

Let $H$ be a bialgebra. If there exists $\lambda \in H^{\ast}$, such that $f\lambda =\varepsilon_{H^{\ast}}(f)\lambda$ for
any $f\in H^{\ast}$, then we call $\lambda$ a left cointegral of $H^{\ast}$. Furthermore, if $\varepsilon_{H^{\ast}}(\lambda)=1$, we call $H$ a cosemisimple bialgebra, and easily see $\lambda^2=\lambda,$ that is, $\lambda$ is idempotent.

So, by the above proposition, we have

\vspace{2mm}

\textbf{Corollary 4.2} Let $H$ be a cosemisimple bialgebra with cointegral $\lambda$, and $M$ a left $H$-comodule. Define two maps
$T: M\rightarrow M$ and $P: H\rightarrow H$ by $T(m)=\lambda(m_{(-1)})m_{(0)}$ and $P(h)=\lambda(h_1)h_2$, for any $m\in M, h\in H$, respectively. Then, $(M,P,T)$ is a Rota-Baxter paired $H$-comodule of weight $-1$.

In particular, $(H, P)$ is a Rota-Baxter coalgbra of weight $-1$ (or by Example 2.2).

 \vspace{2mm}

{\bf 4.2\quad The construction on weak Hopf algebras}

 \vspace{2mm}

\textbf{Definition 4.3} Let $H$ be both an algebra and a coalgebra. Then $H$ is called a {\bf weak bialgebra} in ${\cite{Bohm}}$ if it
satisfies the following conditions:

$(1) \ \Delta(xy)=\Delta(x)\Delta(y),$ for all $x,y\in H$,

$(2)\ \varepsilon(xyz)=\varepsilon(xy_{1})\varepsilon(y_{2}z)=\varepsilon(xy_{2})\varepsilon(y_{1}z),$ for any $x,y,z\in H$,

$(3)\ \Delta^{2}(1_{H})=(\Delta(1_{H})\otimes1_{H})(1_{H}\otimes\Delta(1_{H}))=1_1\otimes 1_21^\prime_1\otimes 1^\prime_2$

\ \ \ \ \ \ \ \ \ \ \ \ \ \ \ \ $=(1_{H}\otimes\Delta(1_{H}))(\Delta(1_{H})\otimes1_{H})=1_1\otimes 1^\prime_11_2\otimes 1^\prime_2,$
\\
where $\Delta(1_{H})=1_{1}\otimes1_{2}=1^\prime_1\otimes 1^\prime_2$ and $\Delta^{2}=(\Delta\otimes id_{H})\circ\Delta$.

Moreover, if there exists a linear map $S : H\rightarrow H$, called
antipode, satisfying the following axioms for all $h\in H$:
$$
h_{1}S(h_{2})=\varepsilon(1_{1}h)1_{2},\ \
S(h_{1})h_{2}=\varepsilon(h1_{2})1_{1},\ \
S(h_{1})h_{2}S(h_{3})=S(h),
$$ then the weak
bialgebra $H$ is called a {\it weak Hopf algebra}.

For any weak bialgebra $H$, defines the maps $\sqcap^L, \sqcap^R: H\rightarrow H$ by the formulas
$$
\sqcap^L(h)=\varepsilon(1_1h)1_2,\ \ \sqcap^R(h)=\varepsilon(h1_2)1_1.
$$

Denote by $H^L$ the image $\sqcap^L$ and by $H^R$ the image $\sqcap^R$, where $H^L$ and $H^R$ are respectively called the target algebra
and the source algebra of the weak bialgebra $H$.

By \cite{Bohm}, if $H$ is a weak Hopf algebra with antipode $S$, we have the following conclusions:

$(W1)$\ \ $\sqcap^L\circ\sqcap^L=\sqcap^L,\ \ \sqcap^R\circ\sqcap^R=\sqcap^R$;

$(W2)$\ \ $\sqcap^L(h_1)\otimes h_2=S(1_1)\otimes 1_2h,\ \ h_1\otimes\sqcap^R(h_2)=h1_1\otimes S(1_2)$, for any $h\in H$;

$(W2^\prime)$\ \ $\sqcap^L(1_1)\otimes 1_2=S(1_1)\otimes 1_2,\ \ 1_1\otimes\sqcap^R(1_2)=1_1\otimes S(1_2)$;

$(W3)$\ \ $\sqcap^L(\sqcap^L(h)g)=\sqcap^L(h)\sqcap^L(g),\ \ \sqcap^R(h\sqcap^R(g))=\sqcap^R(h)\sqcap^R(g)$, for any $h, g\in H.$

Note that $\sqcap^L(1_1)\otimes 1_2$ and  $1_1\otimes \sqcap^R(1_2)$ are separable idempotents of $H^L$ and $H^R$ by Proposition 2.11 in \cite{Bohm}, respectively. So, by $(W2)$, we have

$(W4)$\ \ $h\sqcap^L(1_1)\otimes 1_2=\sqcap^L(1_1)\otimes 1_2h,\ \ h1_1\otimes\sqcap^R(1_2)=1_1\otimes\sqcap^R(1_2)h$, for any $h\in H$.

Again according to \cite{Wang}, we have the following conclusions, that is, for any $x\in H^L, y\in H^R$:

$(W5)$\ \ $\Delta(1)=1_1\otimes 1_2\in H^R\otimes H^L, \ \ xy=yx$;

$(W6)$\ \ $\Delta(x)=1_1x\otimes 1_2,\ \ \Delta(y)=1_1\otimes y1_2$;

$(W7)$\ \ $xS(1_1)\otimes 1_2=S(1_1)\otimes 1_2x,\ \ y1_1\otimes S(1_2)=1_1\otimes S(1_2)y$.

Again by $(W5)$ and $(W6)$, we have

$(W8)$\ \ $\Delta(xy)=1_1x\otimes 1_2y\in H^LH^R\otimes H^LH^R$, for any $x\in H^L, y\in H^R$.

According to $(W8)$, we know that $H^LH^R$ is a subcoalgebra of $H$.

\vspace{2mm}

\textbf{Proposition 4.4} Let $H$ be a weak Hopf algebra with antipode $S$. Then $(H^LH^R, \sqcap^L, \sqcap^L)$ is a Rota-Baxter paired $H$-comodule of weight $-1$, whose comodule structure map is given by the comultiplication $\Delta$ of $H$.

In particular, $(H^LH^R, \sqcap^L)$ is a Rota-Baxter coalgebra of weight $-1$.

{\bf Proof.} By $(W8)$, $H^LH^R$ is a left $H$-comodule via the comultiplication $\Delta$ of $H$. Moreover, for any $x\in H^L, y\in H^R$, we have
\begin{eqnarray*}
(\sqcap^L\o \sqcap^L)\Delta(xy)&\stackrel{(W8)}=&(\sqcap^L\o \sqcap^L)(1_1x\o 1_2y)\\
&=&\sqcap^L(1_1x)\o \sqcap^L(1_2y)\stackrel{(W1)}=\sqcap^L(x1_1)\o \sqcap^L(1_2y)\\
&\stackrel{(W3)}=&x\sqcap^L(1_1)\o 1_2\sqcap^L(y)\\
&\stackrel{(W4)}=&x\sqcap^L(y)\sqcap^L(1_1)\o 1_2,\\
(\sqcap^L\o id)\Delta\sqcap^L(xy)&=&(\sqcap^L\o id)\Delta(x\sqcap^L(y))\\
&=&(\sqcap^L\o id)(x_1\sqcap^L(y)_1\otimes x_2\sqcap^L(y)_2)\\
&\stackrel{(W6)}=&(\sqcap^L\o id)(1_1x1_{1^\prime}\sqcap^L(y)\otimes 1_21_{2^\prime})\\
&=&(\sqcap^L\o id)(1_11_{1^\prime}x\sqcap^L(y)\otimes 1_21_{2^\prime})\\
&=&(\sqcap^L\o id)(1_1x\sqcap^L(y)\otimes 1_2)\\
&=&\sqcap^L(1_1x\sqcap^L(y))\otimes 1_2\\
&\stackrel{(W3)}=&x\sqcap^L(y)\sqcap^L(1_1)\otimes 1_2,\\
(id\o\sqcap^L)\Delta\sqcap^L(xy)&=&(id\o\sqcap^L)\Delta(x\sqcap^L(y))\\
&=&(id\o\sqcap^L)(1_1x\sqcap^L(y)\otimes 1_2)\\
&=&1_1x\sqcap^L(y)\otimes\sqcap^L(1_2)\\
&\stackrel{(W1)}=&1_1x\sqcap^L(y)\otimes 1_2\\
&=&\Delta\sqcap^L(xy).
\end{eqnarray*}

Hence, we get that
\begin{eqnarray*}
(\sqcap^L\o id)\Delta\sqcap^L(xy)+(id\o\sqcap^L)\Delta\sqcap^L(xy)-\Delta\sqcap^L(xy)=(id\o\sqcap^L)\Delta\sqcap^L(xy),
\end{eqnarray*}
that is, $(H^LH^R, \sqcap^L, \sqcap^L)$ is a Rota-Baxter paired $H$-comodule of weight $-1$.
\hfill$\square$

\vspace*{2mm}

In a similar way in Proposition 4.4, we have

\vspace*{2mm}

\textbf{Remark 4.5}  Let $H$ be a weak Hopf algebra with antipode $S$. Then the following hold.

 (1) $(H^LH^R, \sqcap^R, \sqcap^R)$ is a Rota-Baxter paired $H$-comodule of weight $-1$, and so $(H^LH^R, \sqcap^R)$ is a Rota-Baxter coalgebra of weight $-1$;

(2) $(H, \sqcap^L, \sqcap^L)$ is a Rota-Baxter paired $H$-comodule of weight $-1$, and so $(H, \sqcap^L)$ is a Rota-Baxter coalgebras of weight $-1$;

(3) $(H, \sqcap^R, \sqcap^R)$ is a Rota-Baxter paired $H$-comodule of weight $-1$, and so $(H, \sqcap^R)$  is a  Rota-Baxter coalgebra of weight $-1$.

\vspace{2mm}

{\bf 4.3\quad The construction on weak Hopf modules}

 \vspace{2mm}

 In this subsection, we always assume that $H$ is a weak Hopf algebra with antipose $S$. Then, $S$ is both an antimultiplication map and an anticomultiplication map, that is, for any $h, g\in H$,
 $$
 S(hg)=S(g)S(h),\ \ S(1)=1,\ \ \Delta S(h)=S(h_2)\otimes S(h_1),\ \ \varepsilon S(h)=\varepsilon(h),
 $$
 and we have

 $(W9)$\ \ $h_1\otimes \sqcap^L(h_2)=1_1h\otimes 1_2,\ \ \sqcap^R(h_1)\otimes h_2=1_1\otimes h1_2$ for any $h\in H$.

 \vspace{2mm}

\textbf{Definition 4.6} Suppose that $H$ is a weak Hopf algebra with antipode $S$. A {\bf weak right $H$-Hopf module} is a triple $(M,\cdot, \r)$, where $(M,\cdot)$ is a right $H$-module and $(M,\r)$ a right $H$-comodule, such that
$$
\r(m\cdot h)=m_{[0]}\cdot h_1\o m_{[1]}h_2
$$
for any $m\in M, h\in H$.

\vspace*{2mm}

Define a map $T: M\rightarrow M$ given by
\begin{center}
$T(m)=m_{[0]}\cdot S(m_{([1]})$\ \ for\ any\ $m\in M$
\end{center}

Then, according to Proposition 3.8 in \cite{Zheng}, $T$ is idempotent. Moreover, we have

$(W10)$\ \ $\rho T(m)=T(m)\cdot 1_1\otimes 1_2,$ for any $m\in M$.

In fact, for any $m\in M$, we have
\begin{eqnarray*}
\rho T(m)&=&m_{[0][0]}\cdot S(m_{[1]})_1\otimes m_{[0][1]}S(m_{[1]})_2\\
&=&m_{[0]}\cdot S(m_{[1]2})_1\otimes m_{[1]1}S(m_{[1]2})_2\\
&=&m_{[0]}\cdot S(m_{[1]3})\otimes m_{([1]1}S(m_{[1]2})\\
&=&m_{[0]}\cdot S(m_{[1]2})\otimes \sqcap^L(m_{[1]1})\\
&\stackrel{(W2)}=&m_{[0]}\cdot S(1_2m_{[1]})\otimes S(1_1)\\
&=&m_{[0]}\cdot S(m_{[1]})S(1_2)\otimes S(1_1)\\
&=&T(m)\cdot S(1_2)\otimes S(1_1)\\
&=&T(m)\cdot 1_1\otimes 1_2.
\end{eqnarray*}

A weak Hopf algebra $H$ is called a {\bf quantum commutative} if $h_1g\sqcap^R(h_2)=hg$ for any $h,g\in H$. Then, by Proposition 4.1 in \cite{Bagio}, $H$ is quantum commutative if and only if $H^R\subseteq Z(H)$ (the center of $H$).

\vspace{2mm}

\textbf{Proposition 4.7} Let $H$ be a quantum commutative weak Hopf algebra, and $M$ a weak right $H$-Hopf module. Then $(M, \sqcap^L, T)$ is a Rota-Baxter paired $H$-comodule of weight $-1$.

{\bf Proof.} According to $(W9)$, for any $m\in M$, we have
\begin{eqnarray*}
(T\otimes \sqcap^L)\rho (m)&=&T(m_{[0]})\otimes \sqcap^L(m_{[1]})\\
&=&m_{[0][0]}\cdot S(m_{[0][1]})\otimes \sqcap^L(m_{[1]2})\\
&=&m_{[0]}\cdot S(m_{[1]1})\otimes \sqcap^L(m_{[1]2})
\end{eqnarray*}
\begin{eqnarray*}
&\stackrel{(W9)}=&m_{[0]}\cdot S(1_1m_{[1]})\otimes 1_2\\
&=&m_{[0]}\cdot S(m_{[1]})S(1_1)\otimes 1_2\\
&=&T(m)\cdot S(1_1)\otimes 1_2,\\
(id\otimes \sqcap^L)\rho T(m)&\stackrel{(W10)}=&T(m)\cdot 1_1\otimes \sqcap^L(1_2)\\
&\stackrel{(W5)}=&T(m)\cdot 1_1\otimes 1_2\\
&=&\rho T(m),\\
(T\otimes id)\rho T(m)&=&T(T(m)\cdot 1_1)\otimes 1_2\\
&=&T(m)_{[0]}\cdot 1_1S(T(m)_{[1]}1_2)\otimes 1_3\\
&=&T(m)_{[0]}\cdot 1_1S(1_2)S(T(m)_{[1]})\otimes 1_3\\
&=&T(m)_{[0]}\cdot \sqcap^L(1_1)S(T(m)_{[1]})\otimes 1_2\\
&\stackrel{(W2^\prime)}=&T(m)_{[0]}\cdot S(1_1)S(T(m)_{[1]})\otimes 1_2\\
&=&T(m)_{[0]}\cdot S(T(m)_{[1]}1_1)\otimes 1_2\\
&=&T(m)_{[0]}\cdot S(1_1T(m)_{[1]})\otimes 1_2\ \ (H^R\subseteq Z(H))\\
&=&T(m)_{[0]}\cdot S(T(m)_{[1]})S(1_1)\otimes 1_2\\
&=&T(m)\cdot S(1_1)\otimes 1_2,
\end{eqnarray*}
so, we get that
$$
(T\otimes \sqcap^L)\rho (m)=(id\otimes \sqcap^L)\rho T(m)+(T\otimes id)\rho T(m)-\rho T(m).
$$

Hence $(M, T, \sqcap^L)$ is a Rota-Baxter paired $H$-comodule of weight $-1$.

\vspace{2mm}

\textbf{Remark 4.8} Let $H$ be a weak Hopf algebra. Then $H$ is a weak right $H$-Hopf module whose action and coaction are given by its multiplication and comultiplication of $H$. If $H$ is quantum commutative, then, by the above proposition and $h_1S(h_2)=\sqcap^L(h)$ for $h\in H$, we know that $(H, \sqcap^L, \sqcap^L)$ is a Rota-Baxter paired $H$-comodule of weight $-1$.

\vspace{2mm}
{\bf 4.4\quad The construction on dimodules}

 \vspace{2mm}

 In this subsection, we construct Rota-Baxter paired comodules on dimodules.

 \vspace{2mm}

 \textbf{Definition 4.9} Assume that $H$ is a bialgebra. A $k$-module $M$ which is both a left $H$-module and a right $H$-comodule is called a {\bf left, right $H$-dimodule}~\cite{Caenepeel} if for any $h\in H, m\in M$, the following equality holds:
$$
 \rho(h\cdot m)=h\cdot m_{[0]}\otimes m_{[1]},
$$
where $\rho$ is the right $H$-comodule structure map of $M$.

\vspace{2mm}

\textbf{Proposition 4.10} Let $H$ be a bialgebra with an idempotent element $e$, and $M$ a left-right $H$-dimodule. Define a map $T: M\rightarrow M$ by $T(m)=e\c m$. Then $(M, T)$ is a generic Rota-Baxter paired $H$-comodule of weight $-1$.

{\bf Proof.} For any $m\in M$, we have
$$
T^2(m)=T(e\c m)=e^2\c m=e\c m=T(m),
$$
that is, $T$ is idempotent. Since $M$ is a left-right $H$-dimodule, for any $m\in M$, we have
$$
(T\o id)\r(m)=T(m_{[0]})\o m_{[1]}=e\c m_{[0]}\o m_{[1]}=\r(e\c m)=\r T(m),
$$
that is, $T$ is a comodule map. Thus by Theorem 3.1, we know that $(M, T)$ is a generic Rota-Baxter paired $H$-comodule of weight $-1$.
\hfill$\square$

\vspace*{2mm}

\textbf{Remark 4.11} (1) Let $G$ be a finite group. Then, $H=(kG)^\ast=$Hom$_k(kG, k)$ is a Hopf algebra with dual basis $\{p_g|p_g(h)=\delta_{gh}\}$. According to \cite{Cohen}, $p_g$ are orthogonal idempotents, for any $g\in G$. Thus, by the above proposition,  for any left-right $H$-dimodule $M$, $(M, T_g)$ are a generic Rota-Baxter paired $H$-comodule of weight $-1$ for every $g\in G$, where $T_g(m)=p_g\c m$ for $m\in M$.

(2) Let $H$ be a bialgebra. If there is an element $x\in H$ such that $hx=\varepsilon(h)x$ for any $h\in H$, then we call $x$ a {\bf left integral} of $H$.
%

Suppose that $H$ is a finite dimensional semisimple Hopf algebra. Then by \cite[Theorem 5.1.8]{Sweedler}, there exists a non-zero left integral $e$ such that $\varepsilon(e)=1$. It is obvious that $e^2=e$. Hence the following conclusions hold.

(i) If $(H,\mathcal{R})$ is a quasitriangular Hopf algebra. Then, by Example 3.12 in \cite{Zheng}, $(H, \r)$ is a left, right $H$-dimodule, whose action is given by its multiplication and coaction $\r: H\rightarrow H \o H$ given by $\r(h)=h\mathcal{R}_i\o\mathcal{R}_j$. So, by Proposition 4.10, $(H, T)$ is a generic Rota-Baxter paired $H$-comodule of weight $-1$, where $T: H\rightarrow H $ is given by $T(h)=e\c h$.

(ii) Let $M$ be a left, right $H$-dimodule. Define $T(m)=e\cdot m$, for $m\in M$. Then, by Proposition 4.10, $(M, T)$ is a generic Rota-Baxter paired $H$-comodule of weight $-1$.

Note here $T$ is a left $H$-module map: for any $h\in H, m\in M$,
$$
T(h\cdot m)=\varepsilon(h)e\cdot m=he\cdot m=h\cdot T(m).
$$

Again $T$ is idempotent by Proposition 4.10, thus, $(M, T)$ is a generic Rota-Baxter paired $H$-module of weight $-1$ by Theorem 2.4 in \cite{Zheng}. Hence $(M, T)$ is a generic Rota-Baxter paired left, right $H$-dimodule of weight $-1$.

In particular, if $(H,\mathcal{R})$ is a quasitriangular Hopf algebra, then, according to the above conclusions, we know that $(H, T)$ is a generic Rota-Baxter paired left, right $H$-dimodule of weight $-1$.






\vspace{2mm}

{\bf 4.5\quad The construction on relative Hopf modules}

\vspace{2mm}

In this subsection, we construct Rota-Baxter paired comodules on relative Hopf modules.

 \vspace{2mm}

 \textbf{Definition 4.12} Let $H$ be a bialgebra, and $C$ a right
$H$-module coalgebra. A {\bf relative $[C,H]$-Hopf module} in \cite{Zhang} $M$ is a right $C$-comodule which is also a right $H$-module such that the following compatible
condition holds: for all $m \in M$ and $h\in H$,
$$
\rho(m\c h) = m_{[0]}\c h_1 \o m_{[1]}\c h_2.
$$

\textbf{Proposition 4.13}\ \ Let $H$ be a Hopf algebra with an antipode $S$, and $M$ a relative $[C,H]$-Hopf module. If there is a right $H$-module coalgebra map $\phi: C\rightarrow H$, we define $E_C: C\rightarrow C$ and $E_M: M\rightarrow M$ given by
$$
E_C(c)=c_1\cdot S\phi(c_2),
$$
$$
E_M(m)=m_{[0]}\cdot S\phi(m_{[1]}),
$$
for any $c\in C, m\in M$, then, $(M, E_C , E_M )$ is a Rota-Baxter paired right $C$-comodule of weight $-1$.

Here $H$ is a right $H$-module coalgebra whose action is given by its multiplication of $H$.

{\bf Proof.} By the definition of Rota-Baxter paired comodule, we only need to prove that
$$
(E_M\o E_C)\rho(m)=(id\o E_C)\r E_M+(E_M\o id)\r E_M-\r E_M,
$$
for any $m\in M$.

In fact, firstly, for any $m\in M, c\in C$, we have
\begin{eqnarray*}
\r E_M(m)&=&m_{[0][0]}\c (S\phi(m_{[1]}))_1\o m_{[0][1]}\c (S\phi(m_{[1]}))_2\\
&=&m_{[0]}\c (S\phi(m_{[1]2}))_1\o m_{[1]1}\c (S\phi(m_{[1]2}))_2\\
&=&m_{[0]}\c S\phi(m_{[1]3})\o m_{[1]1}\c S\phi(m_{[1]2})\\
&=&m_{[0]}\c S\phi(m_{[1]2})\o E_C(m_{[1]1}),\\
E^2_C(c)&=&E_C(c_1\cdot S\phi(c_2))\\
&=&(c_1\cdot S\phi(c_2))_1\c S\phi ((c_1\cdot S\phi(c_2))_2)\\
&=&c_1\cdot (S\phi(c_4)S\phi (c_2\cdot S\phi(c_3)))\\
&=&c_1\cdot (S\phi(c_4)S(S\phi(c_3))S\phi(c_2))\\
&=&c_1\cdot S(\phi(c_2)S\phi(c_3)\phi(c_4))\\
&=&c_1\cdot S\phi(c_2)=E_C(c).
\end{eqnarray*}

By the above equation, we get
\begin{eqnarray*}
&&(id\o E_C)\r E_M(m)+(E_M\o id)\r E_M(m)-\r E_M(m)\\
&=&(id\o E_C)(m_{[0]}\c S\phi(m_{[1]2})\o E_C(m_{[1]1}))+(E_M\o id)(m_{[0]}\c S\phi(m_{[1]2})\o E_C(m_{[1]1}))\\
&&-m_{[0]}\c S\phi(m_{[1]2})\o E_C(m_{[1]1})\\
&=&m_{[0]}\c S\phi(m_{[1]2})\o E^2_C(m_{[1]1})+E_M(m_{[0]}\c S\phi(m_{[1]2}))\o E_C(m_{[1]1})\\
&&-m_{[0]}\c S\phi(m_{[1]2})\o E_C(m_{[1]1})\\
&=&E_M(m_{[0]}\c S\phi(m_{[1]2}))\o E_C(m_{[1]1}).
\end{eqnarray*}

Finally, we prove that $E_M(m_{[0]}\c S\phi(m_{[1]2}))\o E_C(m_{[1]1})=(E_M\o E_C)\rho(m)$.
\begin{eqnarray*}
&&E_M(m_{[0]}\c S\phi(m_{[1]2}))\o E_C(m_{[1]1})\\
&=&(m_{[0]}\c S\phi(m_{[1]2}))_{[0]}\cdot S\phi((m_{[0]}\c S\phi(m_{[1]2}))_{[1]})\o E_C(m_{[1]1})\\
&=&(m_{[0][0]}\c (S\phi(m_{[1]2}))_1) \cdot S\phi((m_{[0][1]}\c (S\phi(m_{[1]2}))_2)\o E_C(m_{[1]1})\\
&=&m_{[0]}\c (S\phi(m_{[1]4})S(S\phi(m_{[1]3})) S\phi(m_{[1]1})) \o E_C(m_{[1]2})\\
&=&m_{[0]}\c (S(S\phi(m_{[1]3})\phi(m_{[1]4})) S\phi(m_{[1]1})) \o E_C(m_{[1]2})\\
&=&m_{[0]}\c S\phi(m_{[1]1}) \o E_C(m_{[1]2})\\
&=&m_{[0][0]}\c S\phi(m_{[0][1]})\o E_C(m_{[1]})\\
&=&E_M(m_{[0]})\o E_C(m_{[1]}).
\end{eqnarray*}

Hence $(M, E_C , E_M )$ is a Rota-Baxter paired right $C$-comodule of weight $-1$.
\hfill $\square$

\vspace*{2mm}

\textbf{Remark 4.14} Let Let $H$ be a Hopf algebra and and $M$ a relative $[C,H]$-Hopf module, Then, it is easy to show that $H\o M$ is a relative $[C,H]$-Hopf module by
\begin{center}
$\r: H\o M\rightarrow H\o M\o C, ~\r(h\o m)=h_1\o m_{[0]}\o m_{[1]}\c h_2,$
\\
$\cdot: H\o M\o H\rightarrow H\o M,~(h\o m)\c g=hg\o m.$
\end{center}

Then, by Proposition 4.13, $(H\o M, E_C , E_{H\o M})$ is a right Rota-Baxter paired $C$-comodule of weight $-1$ if there exists a right $H$-module coalgebra map $\phi: C\rightarrow H$, where
\begin{center}
$E_C(c)=c_1\cdot S\phi(c_2)$,\\
$E_{H\o M}(h\o m)=h_1S\phi(m_{[1]}\c h_2)\o m_{[0]}.$
\end{center}

\vspace*{2mm}

{\bf 4.6\quad The construction on Rota-Baxter paired comodules}

 \vspace{2mm}

 \textbf{Proposition 4.15}\ \  Let $(M,P,T)$ be a Rota-Baxter paired $C$-comodule of weight $\lambda$. Define
$$\overline{P}=-P-\lambda id, ~~\overline{T}=-T-\lambda id $$
Then $(M,\overline{P},\overline{T})$ is also a Rota-Baxter paired $C$-comodule of weight $\lambda$.

{\bf Proof.} We have only to verify that
$$
(\overline{P}\otimes \overline{T})\rho=(\overline{P}\otimes id)\rho \overline{T}+(id\otimes \overline{T})\rho\overline{T}+\lambda\rho \overline{T}.
$$

Actually, for any $m\in M$, we have
\begin{eqnarray*}
&&(\overline{P}\otimes \overline{T})\rho(m)=(\overline{P}\otimes \overline{T})(m_{(-1)}\otimes m_{(0)})=\overline{P}(m_{(-1)})\otimes \overline{T}((m_{(0)}))\\
&=&(-P(m_{(-1)})-\lambda m_{(-1)})\otimes(-T(m_{(0)})-\lambda m_{(0)})\\
&=&P(m_{(-1)})\otimes T(m_{(0)})+\lambda P(m_{(-1)})\otimes m_{(0)}+\lambda m_{(-1)}\otimes T(m_{(0)})+\lambda^2m_{(-1)}\otimes m_{(0)}\\
&=&P(T(m)_{(-1)})\otimes T(m)_{(0)}+T(m)_{(-1)}\otimes T(T(m)_{(0)})+\lambda T(m)_{(-1)}\otimes T(m)_{(0)}\\
&&+\lambda P(m_{(-1)})\otimes m_{(0)}+\lambda m_{(-1)}\otimes T(m_{(0)})+\lambda^2m_{(-1)}\otimes m_{(0)},\\
&&((\overline{P}\otimes id)\rho \overline{T}+(id\otimes \overline{T})\rho \overline{T}+\lambda\rho \overline{T})(m)=(\overline{P}\otimes id+id\otimes \overline{T}+\lambda)\rho(-T(m)-\lambda m)\\
&=&-(\overline{P}\otimes id)(T(m)_{(-1)}\otimes T(m)_{(0)})-\lambda(\overline{P}\otimes id)(m_{(-1)}\otimes m_{(0)})\\
&&-(id\otimes \overline{T})(T(m)_{(-1)}\otimes T(m)_{(0)})-\lambda(id\otimes \overline{T})(m_{(-1)}\otimes m_{(0)})\\
&&-\lambda T(m)_{(-1)}\otimes T(m)_{(0)}-\lambda^2m_{(-1)}\otimes m_{(0)}\\
&=&-\overline{P}(T(m)_{(-1)})\otimes T(m)_{(0)}-\lambda \overline{P}(m_{(-1)})\otimes m_{(0)}-T(m)_{(-1)}\otimes \overline{T}(T(m)_{(0)})\\
&&-\lambda m_{(-1)}\otimes \overline{T}(m_{(0)})-\lambda T(m)_{(-1)}\otimes T(m)_{(0)}-\lambda^2m_{(-1)}\otimes m_{(0)}\\
&=&(P(T(m)_{(-1)})+\lambda T(m)_{(-1)})\otimes T(m)_{(0)}-\lambda(-P(m_{(-1)})-\lambda m_{(-1)})\otimes m_{(0)}\\
&&-T(m)_{(-1)}\otimes(-T(T(m)_{(0)})-\lambda T(m)_{(0)})-\lambda m_{(-1)}\otimes(-T(m_{(0)})-\lambda m_{(0)})\\
&&-\lambda T(m)_{(-1)}\otimes T(m)_{(0)}-\lambda^2m_{(-1)}\otimes m_{(0)}\\
&=&P(T(m)_{(-1)})\otimes T(m)_{(0)}+\lambda T(m)_{(-1)}\otimes T(m)_{(0)}+\lambda P(m_{(-1)})\otimes m_{(0)}\\
&&+T(m)_{(-1)}\otimes T(T(m)_{(0)})+\lambda m_{(-1)}\otimes T(m_{(0)})+\lambda^2m_{(-1)}\otimes m_{(0)}
\end{eqnarray*}
as desired.
\hfill $\square$

\vspace*{2mm}

\textbf{Proposition 4.16}\ \ Let $(C,P)$ be a Rota-Baxter coalgebra of weight $\lambda$, and $(M,P,T)$ a Rota-Baxter paired comodule of weight $\lambda$. Define another comultiplication $\Delta^\prime$ on $C$ by
$$\Delta^\prime=(id\o P)\Delta+(P\o id)\Delta+\lambda \Delta,
$$
and another operation of $M$ by $$\r^\prime=(P\o id)\r+(id\o T)\r+\lambda \r.
$$

Then the following conclusion hold.

(a) $(C,\Delta^\prime,P)$ is also a (noncounitary) Rota-Baxter coalgebra of weight $\lambda$.

(b) $\r^\prime T=(P\o T)\r$.

(c) $(M,\r^\prime)$ is a noncounitary $(C,\Delta^\prime)$-comodule.

(d) $(M,P,T)$ is a Rota-Baxter paired $(C,\Delta^\prime)$-comodule of weight $\lambda$, whose comodule structure map is given by $\r^\prime$.

{\bf Proof.} (a) Since $(C,\Delta,P)$ is a Rota-Baxter coalgebra of weight $\lambda$, we have
$$
(P\otimes P)\Delta=(P\otimes id)\Delta P+(id\otimes P)\Delta P+\lambda \Delta P.
$$

Hence we get $\Delta^\prime P=(P\otimes P)\Delta$.

Moreover, for any $c\in C$, we have
\begin{eqnarray*}
(P\otimes P)\Delta^\prime(c)&=&(P\otimes P)((P\otimes id)\Delta+(id\otimes P)\Delta+\lambda\Delta)(c)\\
&=&(P\otimes P)(P(c_{1})\otimes c_{2}+c_{1}\otimes P(c_{2})+\lambda c_{1}\otimes c_{2})\\
&=&P^2(c_{1})\otimes P(c_{2})+P(c_{1})\otimes P^2(c_{2})+\lambda P(c_{1})\otimes P(c_{2})\\
&=&(P\otimes id+id\otimes P+\lambda)(P(c_{1})\otimes P(c_{2}))\\
&=&(P\otimes id+id\otimes P+\lambda)(P\otimes P)\Delta(c)\\
&=&((P\otimes id)\Delta^\prime P+(id \otimes P)\Delta^\prime P+\lambda\Delta^\prime P)(c).
\end{eqnarray*}

(b) It is direct to check $\r^\prime T=(P\o id)\r T+(id\o T)\r T+\lambda \r T=(P\o T)\r$.

(c) We only need to prove that $(id\o \r^\prime)\r^\prime=(\D^\prime\o id)\r^\prime$: by $\D^\prime P=(P\o P)\D$ and $(b)$, for any $m\in M$, we have
\begin{eqnarray*}
&&(id\o \r^\prime)\r^\prime(m)=P(m_{(-1)})\o \r^\prime(m_{(0)})+m_{(-1)}\o \r^\prime(T(m_{(0)}))+\lambda m_{(-1)}\o \r^\prime(m_{(0)})\\
&=&P(m_{(-1)})\o P(m_{(0)(-1)})\o m_{(0)(0)}+P(m_{(-1)})\o m_{(0)(-1)}\o T(m_{(0)(0)})\\
&&+P(m_{(-1)})\o \lambda m_{(0)(-1)}\o m_{(0)(0)}+m_{(-1)}\o \r^\prime(T(m_{(0)}))+\lambda m_{(-1)}\o P(m_{(0)(-1)})\o m_{(0)(0)}\\
&&+\lambda m_{(-1)}\o m_{(0)(-1)}\o T(m_{(0)(0)})+\lambda m_{(-1)}\o \lambda m_{(0)(-1)}\o m_{(0)(0)}\\
&=&P(m_{(-1)1})\o P(m_{(-1)2})\o m_{(0)}+P(m_{(-1)1})\o m_{(-1)2}\o T(m_{(0)})\\
&&\underbrace{+P(m_{(-1)1})\o \lambda m_{(-1)2}\o m_{(0)}}+m_{(-1)}\o \r^\prime(T(m_{(0)}))\underbrace{+\lambda m_{(-1)1}\o P(m_{(-1)2})\o m_{(0)}}\\
&&+\lambda m_{(-1)1}\o m_{(-1)2}\o T(m_{(0)})+\underbrace{\lambda m_{(-1)1}\o \lambda m_{(-1)2}\o m_{(0)}}\\
&=&P(m_{(-1)1})\o P(m_{(-1)2})\o m_{(0)}+P(m_{(-1)1})\o m_{(-1)2}\o T(m_{(0)})+\underbrace{\D^\prime(m_{(-1)})\o \lambda m_{(0)}}\\
&&+m_{(-1)}\o \r^\prime(T(m_{(0)}))+\lambda m_{(-1)1}\o m_{(-1)2}\o T(m_{(0)})
\end{eqnarray*}
\begin{eqnarray*}
&=&P(m_{(-1)1})\o P(m_{(-1)2})\o m_{(0)}+P(m_{(-1)1})\o m_{(-1)2}\o T(m_{(0)})+\D^\prime(m_{(-1)})\o \lambda m_{(0)}\\
&&+m_{(-1)}\o (P\o T)\r(m_{(0)}))+\lambda m_{(-1)1}\o m_{(-1)2}\o T(m_{(0)})\\
&=&P(m_{(-1)1})\o P(m_{(-1)2})\o m_{(0)}\underbrace{+P(m_{(-1)1})\o m_{(-1)2}\o T(m_{(0)})}+\D^\prime(m_{(-1)})\o \lambda m_{(0)}\\
&&+\underbrace{m_{(-1)1}\o P(m_{(-1)2})\o T(m_{(0)})}+\underbrace{\lambda m_{(-1)1}\o m_{(-1)2}\o T(m_{(0)})}\\
&=&P(m_{(-1)1})\o P(m_{(-1)2})\o m_{(0)}+\underbrace{\D^\prime(m_{(-1)})\o T(m_{(0)})}+\D^\prime(m_{(-1)})\o \lambda m_{(0)}\\
&=&\D^\prime P(m_{(-1)})\o m_{(0)}+\D^\prime(m_{(-1)})\o T(m_{(0)})+\D^\prime(m_{(-1)})\o \lambda m_{(0)}\\
&=&(\D^\prime\o id)((P\o id)\r(m)+(id\o T)\r(m)+\lambda \r(m))\\
&=&(\D^\prime\o id)\r^\prime(m).
\end{eqnarray*}

(d) By (b) and (c), we have only to prove that
\begin{eqnarray*}
&&(P\o T)\r^\prime(m)=(P\o T)(P(m_{(-1)})\o m_{(0)}+m_{(-1)}\o T(m_{(0)})+\lambda m_{(-1)}\o m_{(0)})\\
&=&P^2(m_{(-1)})\o T(m_{(0)})+P(m_{(-1)})\o T^2(m_{(0)})+\lambda P(m_{(-1)})\o T(m_{(0)})\\
&=&(P\o id+id\o T+\lambda)(P\o T)\r(m)\\
&=&(P\o id+id\o T+\lambda)\r^\prime T(m)\\
&=&(P\o id)\r^\prime T(m)+(id\o T)\r^\prime T(m)+\lambda\r^\prime T(m)
\end{eqnarray*}
 for any $m\in M$, so, (d) holds.\hfill $\square$

\vspace*{2mm}

By the above propositions, we easily get the following corollary.

\vspace*{2mm}

\textbf{Corollary 4.17}\ \  Let $(M,P,T)$ be a Rota-Baxter paired $C$-comodule of weight $\lambda$. Then $(M,\overline{P},\overline{T})$ is also a Rota-Baxter paired $(C,\Delta^\prime)$-comodule of weight $\lambda$.

 Here the coaction $\r^\prime$ of $M$ and the comultiplication $\Delta^\prime$ of $C$ are defined in Proposition 4.16, and $\overline{P}, \overline{T}$ are defined in Proposition 4.12.

\vspace*{2mm}

\begin{center}
{\bf \S5\quad From Rota-Baxter paired comodules to Pre-Lie comodules}
\end{center}

In this section, we mainly construct pre-Lie comodules from Rota-Baxter paired comodules.
\vspace*{2mm}

\textbf{Definition 5.1} A {\bf pre-Lie coalgebra} is $(C,\D)$ consisting of a linear space $C$, a linear map $\D : C \rightarrow C \o C$ and satisfying
$$
\D_{C}-\Phi_{(12)}\D_{C}=0,
$$
where $\D_{C}=(\D\o id)\D-(id\o \D)\D$ and $\Phi_{(12)}(c_1\o c_2\o c_3)=c_2\o c_1\o c_3$.
\vspace*{2mm}

\textbf{Definition 5.2} Let $(C,\D)$ be a pre-Lie coalgebra. A left {\bf $C$-pre-Lie comodule} $(M, \rho)$ is a space $M$ together with a map $\r: M\rightarrow C\o M$, such that
$$
\r_{M}-(\tau\o id)\r_M=0,
$$
where $\r_M=(id\o \r)\r-(\D\o id)\r$, and $\tau(c\o d)=d\o c$, for any $c,d\in C$.
\vspace*{2mm}

\textbf{Lemma 5.3} Let $(C, Q)$ be a Rota-Baxter coalgebra of weight $-1$. Define the
operation $\widetilde{\D}$ on $C$ by
$$
\widetilde{\D}(c)=Q(c_1)\o c_2-Q(c_2)\o c_1-c_1\o c_2.
$$

Then, by \cite{Ma}, $\widetilde{C}=(C,\widetilde{\D})$ is a pre-Lie coalgebra.

\textbf{Proposition 5.4} Let $(C, Q)$ be a Rota-Baxter coalgebra of weight $-1$, and $(M, Q, T)$ a Rota-Baxter paired $C$-comodule of weight $-1$.
Define a map $\widetilde{\r}: M\rightarrow C\o M$ by
$$\widetilde{\r}(m)=Q(m_{(-1)})\o m_{(0)}+m_{(-1)}\o T(m_{(0)})-m_{(-1)}\o m_{(0)}.
$$

Then $(M,\widetilde{\r})$ is a left $\widetilde{C}$-pre-Lie comodule, where $\widetilde{C}$ is defined as in Lemma 5.3.

{\bf Proof.} By Lemma 5.3, we know $\widetilde{C}=(C,\widetilde{\D})$ is a pre-Lie coalgebra, so, we only need to prove $\widetilde{\r}_M-(\tau\o id)\widetilde{\r}_M=0$ holds.

As a matter of fact, for any $m\in M$, we have
\begin{eqnarray*}
\widetilde{\r}_M(m)&=&(id\o \widetilde{\r})\widetilde{\r}(m)-(\widetilde{\D}\o id)\widetilde{\r}(m)\\
&=&(id\o \widetilde{\r})(Q(m_{(-1)})\o m_{(0)}+m_{(-1)}\o T(m_{(0)})-m_{(-1)}\o m_{(0)})\\
&&-(\widetilde{\D}\o id)(Q(m_{(-1)})\o m_{(0)}+m_{(-1)}\o T(m_{(0)})-m_{(-1)}\o m_{(0)})\\
&=& Q(m_{(-1)})\o \widetilde{\r}(m_{(0)})+m_{(-1)}\o \widetilde{\r}(T(m_{(0)}))-m_{(-1)}\o \widetilde{\r}(m_{(0)})\\
&&-\widetilde{\D}(Q(m_{(-1)}))\o m_{(0)}-\widetilde{\D}(m_{(-1)})\o T(m_{(0)})+\widetilde{\D}(m_{(-1)})\o m_{(0)}\\
&=& Q(m_{(-1)})\o Q(m_{(0)(-1)})\o m_{(0)(0)}+Q(m_{(-1)})\o m_{(0)(-1)}\o T(m_{(0)(0)})\\
&&-Q(m_{(-1)})\o m_{(0)(-1)}\o m_{(0)(0)}+m_{(-1)}\o Q(T(m_{(0)})_{(-1)})\o T(m_{(0)})_{(0)}\\
&&+m_{(-1)}\o T(m_{(0)})_{(-1)}\o T(T(m_{(0)})_{(0)})-m_{(-1)}\o T(m_{(0)})_{(-1)}\o T(m_{(0)})_{(0)}\\
&&-m_{(-1)}\o Q(m_{(0)(-1)})\o m_{(0)(0)}-m_{(-1)}\o m_{(0)(-1)}\o T(m_{(0)(0)})\\
&&+m_{(-1)}\o m_{(0)(-1)}\o m_{(0)(0)}-Q(Q(m_{(-1)})_{1})\o Q(m_{(-1)})_{2}\o m_{(0)}\\
&&+Q(Q(m_{(-1)})_{2})\o Q(m_{(-1)})_{1}\o m_{(0)}+Q(m_{(-1)})_{1}\o Q(m_{(-1)})_{2}\o m_{(0)}\\
&&-Q(m_{(-1)1})\o m_{(-1)2}\o T(m_{(0)})+Q(m_{(-1)2})\o m_{(-1)1}\o T(m_{(0)})\\
&&+m_{(-1)1}\o m_{(-1)2}\o T(m_{(0)})+Q(m_{(-1)1})\o m_{(-1)2}\o m_{(0)}\\
&&-Q(m_{(-1)2})\o m_{(-1)1}\o m_{(0)}-m_{(-1)1}\o m_{(-1)2}\o m_{(0)}\\
&=& Q(m_{(-1)1})\o Q(m_{(-1)2})\o m_{(0)}+Q(m_{(-1)1})\o m_{(-1)2}\o T(m_{(0)})\\
&&-Q(m_{(-1)1})\o m_{(-1)2}\o m_{((0)}+m_{(-1)}\o Q(T(m_{(0)})_{(-1)})\o T(m_{(0)})_{(0)}\\
&&+m_{(-1)}\o T(m_{(0)})_{(-1)}\o T(T(m_{(0)})_{(0)})-m_{(-1)}\o T(m_{(0)})_{(-1)}\o T(m_{(0)})_{(0)}\\
&&-m_{(-1)}\o Q(m_{(0)(-1)})\o m_{(0)(0)}-m_{(-1)1}\o m_{(-1)2}\o T(m_{(0)})\\
&&+m_{(-1)1}\o m_{(-1)2}\o m_{(0)}-Q(Q(m_{(-1)})_{1})\o Q(m_{(-1)})_{2}\o m_{(0)}\\
&&+Q(Q(m_{(-1)})_{2})\o Q(m_{(-1)})_{1}\o m_{(0)}+Q(m_{(-1)})_{1}\o Q(m_{(-1)})_{2}\o m_{(0)}\\
&&-Q(m_{(-1)1})\o m_{(-1)2}\o T(m_{(0)})+Q(m_{(-1)2})\o m_{(-1)1}\o T(m_{(0)})\\
&&+m_{(-1)1}\o m_{(-1)2}\o T(m_{(0)})+Q(m_{(-1)1})\o m_{(-1)2}\o m_{(0)}\\
&&-Q(m_{(-1)2})\o m_{(-1)1}\o m_{(0)}-m_{(-1)1}\o m_{(-1)2}\o m_{(0)}\\
&=& Q(m_{(-1)1})\o Q(m_{(-1)2})\o m_{(0)}+m_{(-1)}\o Q(T(m_{(0)})_{(-1)})\o T(m_{(0)})_{(0)}\\
&&+m_{(-1)}\o T(m_{(0)})_{(-1)}\o T(T(m_{(0)})_{(0)})-m_{(-1)}\o T(m_{(0)})_{(-1)}\o T(m_{(0)})_{(0)}\\
&&-m_{(-1)}\o Q(m_{(0)(-1)})\o m_{(0)(0)}-Q(Q(m_{(-1)})_{1})\o Q(m_{(-1)})_{2}\o m_{(0)}\\
\end{eqnarray*}
\begin{eqnarray*}
&&+Q(Q(m_{(-1)})_{2})\o Q(m_{(-1)})_{1}\o m_{(0)}+Q(m_{(-1)})_{1}\o Q(m_{(-1)})_{2}\o m_{(0)}\\
&&+Q(m_{(0)(-1)})\o m_{(-1)}\o T(m_{(0)(0)})-Q(m_{(0)(-1)})\o m_{(-1)}\o m_{(0)(0)}.
\end{eqnarray*}

So, by the above equality and $(C, Q)$ being a Rota-Baxter coalgebra of weight $-1$, we easily prove that $\widetilde{\r}_M(m)=(\tau\o id)\widetilde{\r}_M(m)$ for any $m\in M$. Hence $(M,\widetilde{\r})$ is a left $\widetilde{C}$-pre-Lie comodule.
\hfill $\square$

\vspace*{2mm}

\textbf{Lemma 5.5} Let $(C,\D,Q)$ be a Rota-Baxter coalgebra of weight $0$. Define the
operation $\widetilde{\D}$ on $C$ by
$$
\widetilde{\D}(c)=Q(c_1)\o c_2-Q(c_2)\o c_1.
$$

Then, by \cite{Ma}, $\widetilde{C}=(C,\widetilde{\D})$ is a pre-Lie coalgebra.

\vspace*{2mm}

According to Lemma 5.5, we can prove the following proposition using a similar way as in Proposition 5.4.

\vspace*{2mm}

\textbf{Proposition 5.6} Let $(C, Q)$ be a Rota-Baxter coalgebra of weight $0$, and $(M, Q, T)$ a Rota-Baxter paired $C$-comodule of weight $0$. Define a map $\widetilde{\r}: M\rightarrow C\o M$ by
$$\widetilde{\r}(m)=Q(m_{(-1)})\o m_{(0)}+m_{(-1)}\o T(m_{(0)})
,$$
 for any $m\in M$.
Then $(M,\widetilde{\r})$ is a left $\widetilde{C}$-pre-Lie comodule, where $\widetilde{C}$ is defined as in Lemma 5.5.

\vspace*{2mm}

\begin{center}
{\bf \S6\quad Rota-Baxter paired Hopf module}
\end{center}

In this section, we will combine Rota-Baxter paired modules and Rota-Baxter paired comodules, and introduce the conception of Rota-Baxter pared Hopf modules, and give the structure theorem of Rota-Baxter pared Hopf modules.

\vspace*{2mm}

\textbf{Definition 6.1} Let $H$ be a bialgera, and $M$ a left $H$-Hopf module. A triple $(M, P, T)$ is called a {\bf Rota-Baxter paired left $H$-Hopf module of weight $\lambda$}, if $(M, P, T)$ is both a Rota-Baxter paired left $H$-module of weight $\lambda$, and a Rota-Baxter paired left $H$-comodule of weight $\lambda$.

\vspace*{2mm}

A {\bf Rota-Baxter $H$-Hopf submodule} $N$ of a Rota-Baxter paired $H$-Hopf-module $(M, P, T)$ is an
$H$-Hopf submodule of $M$ such that $T(N)\subseteq N$. Then $(N,P,T)$ is a Rota-Baxter paired $H$-Hopf module.

Let $(M, P, T)$ and $(M^\prime, P^\prime, T^\prime)$ be Rota-Baxter paired $H$-Hopf modules of the same weight $\lambda$. A {\bf Rota-Baxter $H$-Hopf module map} $f: (M, P, T)\rightarrow (M^\prime , P^\prime, T^\prime)$ of weight $\lambda$ is a Hopf module map such that $f\circ T=T^\prime\circ f$.

\vspace*{2mm}

\textbf{Example 6.2} (1) Let $H$ be a bialgebra. Then, $H$ is not only an augmented coalgebra (there is a coalgebra map $\mu: k\rightarrow H$) and an augmented algebra (there is an algebra map $\varepsilon: H\rightarrow k$).

Define a map $P: H\rightarrow H$ given by $P(h)=\varepsilon(h)1_H$. Then, by Example 2.2, $(H, P)$ is a Rota-Baxter coalgebra of weight $-1$, and a Rota-Baxter algebra of weight $-1$ by Example 2.1 in \cite{Zheng}. So, $(H, P, P)$ is not only a Rota-Baxter paired $H$-comodule of weight $-1$, and a Rota-Baxter paired $H$-module of weight $-1$.

It is obvious that $H$ is a right $H$-Hopf module via its multiplication and its comultiplication. Hence $(H, P, P)$ is a Rota-Baxter paired $H$-Hopf module of weight $-1$.

(2) Let $H$ be a bialgebra, and $(H, P, P)$ a Rota-Baxter bialgebra of weight $(\lambda, \lambda)$ given in \cite{Ma}, that is, $(H, P)$ is both a Rota-Baxter algebra of weight $\lambda$, and $(H, P)$ a Rota-Baxter coalgebra of weight $\lambda$. Then, $(H, P, P)$ is both a Rota-Baxter paired $H$-module of weight $\lambda$, whose action is given by the multiplication  of $H$, and $(H, P, P)$ a Rota-Baxter paired $H$-comodule of weight $\lambda$, whose coaction is given by the comultiplication of $H$. So, $(H, P, P)$ is  a Rota-Baxter paired $H$-Hopf module of weight $\lambda$.

(3) Let $H$ be a quantum commutative weak Hopf algebra, and $M$ a weak $H$-Hopf module. Then, by Remark 3.17 in \cite{Zheng}, $(M, \sqcap^L, T)$ is a Rota-Baxter paired $H$-module of weight $-1$, where $T$ is given in Proposition 4.7.

Again by Proposition 4.7, $(M, \sqcap^L, T)$ is also a Rota-Baxter paired $H$-comodule of weight $-1$. So, $(M, \sqcap^L, T)$ is a Rota-Baxter paired $H$-Hopf module of weight $-1$.

In particular, $(H, \sqcap^L, \sqcap^L)$ is a Rota-Baxter paired $H$-Hopf module of weight $-1$ for every quantum commutative weak Hopf algebra $H$.

\vspace*{2mm}

Combining Theorem 3.1 and Theorem 2.4 in \cite{Zheng}, we get the following result.

\vspace*{2mm}

\textbf{Proposition 6.3} Let $H$ be a bialgebra, and $M$ a left $H$-Hopf module. Suppose that there is a Hopf module map $T$ from $M$ to $M$. Then the following are equivalent:

(1) $(M,T)$ is a generic Rota-Baxter paired $H$-Hopf module of weight $\lambda$;

(2) There is a linear operator $P:H\rightarrow H$ such that $(M,P,T)$ is a Rota-Baxter paired $H$-Hopf module of weight $\lambda$;

(3) $T$ is quasi-idempotent of weight $\lambda$.

\vspace*{2mm}

\textbf{Example 6.4} Let $H$ be a bialgebra, and $C$ a coalgebra. Then $H\o C$ is left $H$-Hopf module by $h\c (g\o c)=hg\o c$ and $\r(h\o c)=h_1\o h_2\o c$, for any $h,g\in H, c\in C$. Define a map $T: H\o C\rightarrow H\o C$ by $T(h\o c)=h\o \v(c)e$, where $e\in C$ satisfying $\v(e)=1$. It isn't difficult to prove that $T^2=T$ and $T$ a Hopf module map, so by Proposition 6.3, $(H\o C,T)$ is a generic Rota-Baxter paired $H$-Hopf module of weight $-1$.

\vspace*{2mm}

In what follows, we give the structure theorem of the generic Rota-Baxter  paired Hopf module.

\vspace*{2mm}

\textbf{Theorem 6.5} Let $H$ be a Hopf algebra, and $(M,T)$ a generic Rota-Baxter pared $H$-Hopf module of weight $\lambda$  in Proposition 6.3. Then, there is an isomorphism:
$$
(M, T)\cong (H\o M^{coH}, T^\prime)
$$
as generic Rota-Baxter left $H$-Hopf modules of weight $\lambda$, where $T^\prime$ is defined by
$$T^\prime(h\o m)=h\o T(m), \ \ h\in H, m\in M^{coH},
$$
 and  $M^{coH}=\{m\in M\mid \r(m)=1\o m\}$, and $H\o M^{coH}$ is left $H$-Hopf module by $h\c (g\o m)=hg\o m$ and $\r(h\o m)=h_1\o h_2\o m$, for any $h,g\in H, m\in M^{coH}$.

{\bf Proof.} Since $T$ is a left $H$-comodule map, we easily see that $T(m)\in M^{coH}$, for $m\in M^{coH}$. Hence $T^\prime$ is well defined.

According to \cite{Sweedler}, it is obvious that $H\o M^{coH}$ is a left $H$-Hopf module.

In what follows, by Proposition 6.3, we prove that $(H\o M^{coH}, T')$ is a generic Rota-Baxter left $H$-Hopf module of weight $\lambda$.

As a matter of fact, for any $h\in H, m\in M^{coH}$, we have
$${T^\prime}^2(h\o m)=T^\prime(h\o T(m))=h\o T^2(m)=-\lambda h\o T(m)=-\lambda T^\prime(h\o m),
$$
so $T^\prime$ is quasi-idempotent of weight $\lambda$. Again, for any $h, g\in H, m\in M^{coH}$, we get
\begin{eqnarray*}
T^\prime(h\c (g\o m))&=&T^\prime(hg\o m)=hg\o T(m)\\
&=&h\c T^\prime(g\o m),\\
T^\prime(h\o m)_{(-1)}\o T^\prime(h\o m)_{(0)}&=&(h\o T(m))_{(-1)}\o (h\o T(m))_{(0)}\\
&=&h_1\o (h_2\o T(m))=h_1\o T^\prime(h_2\o m)\\
&=&(h\o m)_{(-1)}\o T^\prime((h\o m)_{(0)}),
\end{eqnarray*}
so, $T^\prime$ is a left $H$-Hopf module map. Hence $(H\o M^{coH}, T^\prime)$ is a generic Rota-Baxter left $H$-Hopf module of weight $\lambda$ by Proposition 6.3.

According to Theorem 4.1.1 in \cite{Sweedler}, we have an $H$-Hopf module isomorphisms as follows:
$$
\alpha: H\o M^{coH}\rightarrow M, ~~h\o m\mapsto h\cdot m,
$$
with the inverse
$$
\beta: M\rightarrow H\o M^{coH}, ~~m\mapsto m_{(-1)}\o E_M(m_{(0)}),
$$
where $E_M(m)$ is given by $S(m_{(-1)})\c m_{(0)}$, for $m\in M$.

Moreover, for any $m\in M$, we obtain
\begin{eqnarray*}
T^\prime\circ \beta(m)&=&T^\prime(m_{(-1)}\o E_M(m_{(0)}))\\
&=&m_{(-1)}\o T(E_M(m_{(0)}))\\
&=&m_{(-1)}\o T(S(m_{(0)(-1)})\c m_{(0)(0)})\\
&=&m_{(-1)}\o S(m_{(0)(-1)})\c T(m_{(0)(0)})\\
&=&m_{(-1)}\o S(T(m_{(0)})_{(-1)})\c T(m_{(0)})_{(0)}\\
&=&m_{(-1)}\o E_M(T(m_{(0)}))\\
&=&T(m)_{(-1)}\o E_M(T(m)_{(0)}))\\
&=&\beta\circ T(m),
\end{eqnarray*}
so $T^\prime\circ \beta=\beta\circ T$. In a similar way, we can prove that $\alpha\circ T^\prime=T\circ \alpha$. Hence $(M, T)\cong (H\o M^{coH}, T^\prime)$ as generic Rota-Baxter left $H$-Hopf module of weight $\lambda$. \hfill $\square$

\vspace{5mm}

\makeatletter
\newcommand{\adjustmybblparameters}{\setlength{\itemsep}{0\baselineskip}\setlength{\parsep}{0pt}}
\let\ORIGINALlatex@openbib@code=\@openbib@code
\renewcommand{\@openbib@code}{\ORIGINALlatex@openbib@code\adjustmybblparameters}
\makeatother

\end{document}